\documentclass{article}
\usepackage[english]{babel}
\usepackage{amsmath}
\usepackage{amsthm}
\usepackage{amssymb}
\usepackage{amscd}
\usepackage[latin1]{inputenc}

\newtheorem{teor}{Theorem}[section]
\newtheorem{defi}{Definition}
\newtheorem{lema}[teor]{Lemma}
\newtheorem{prop}[teor]{Proposition}
\newtheorem{cor}[teor]{Corollary}
\newtheorem{rem}[teor]{Remark}
\newtheorem{ejem}[teor]{Example}
\newtheorem{ejems}[teor]{Examples}
\newtheorem{ques}[teor]{Question}

\begin{document}

\title{Killing of supports  on graded  algebras}

\author{Roberto Martínez Villa\\
Instituto de Matemáticas, UNAM, AP 61-3 \\
58089 Morelia, Michoacán\\ MEXICO \\ {\it mvilla@matmor.unam.mx} \vspace*{0.5cm} \and Manuel Saorín\\ Departamento de Matemáticas\\
Universidad de Murcia, Aptdo. 4021\\
30100 Espinardo, Murcia\\
SPAIN\\ {\it msaorinc@um.es}}

\date{}

\thanks{
The first named author thanks CONACYT for funding the research
project. The second one thanks the D.G.I. of the Spanish Ministry
of Education and Science and the Fundación "Séneca" of Murcia for
their financial support}

\maketitle

\begin{abstract}

{\bf Killing of supports along subsets $\mathcal{U}$ of a group
$G$ and regradings along certain maps of groups $\varphi
:G'\longrightarrow G$ are studied, in the context of group-graded
algebras. We show that, under precise conditions on $\mathcal{U}$
and $\varphi$, the module theories over the initial and the final
algebras are functorially well-connected. Special attention is
paid to  $G=\mathbf{Z}$, in which case the results can be applied
to $n$-Koszul algebras
 }

\end{abstract}

\section{Introduction}

In a recent paper (\cite{GMMVZ}, see also \cite{BM}),   the
authors showed that if $\Lambda$ is a $n$-Koszul graded algebra
and $\Lambda^!$ is its $n$-homogeneous dual, then the Yoneda
algebra $E=E(\Lambda )$ of $\Lambda$ is obtained as follows. One
takes $\mathcal{U}=n\mathbf{Z}\cup (n\mathbf{Z}+1)$,  kills the
part of the support of $\Lambda^!$ in
$\mathbf{Z}\setminus\mathcal{U}$ and then obtains a new
$\mathbf{Z}$-graded algebra $\Lambda^!_\mathcal{U}$ by keeping the
same multiplication of homogeneous elements,   when its degree
falls into $\mathcal{U}$, and making it zero otherwise.   Then $E$
is obtained from $\Lambda^!_\mathcal{U}$ by regrading along
certain function $\delta :\mathbf{Z}\longrightarrow\mathbf{Z}$
which has image $Im(\delta )=\mathcal{U}$.

Motivated by the above fact, we consider an arbitrary subset
$\mathcal{U}$ of a group $G$ and define a multiplication of
nonzero  homogeneous elements  on
$A_\mathcal{U}=\oplus_{u\in\mathcal{U}}A_u$ by the rule $a\cdot
b=ab$, in case $deg(a)deg(b)\in\mathcal{U}$, and $a\cdot b=0$
otherwise. We then extend by linearity and  call this operation
the {\bf support-restricted multiplication} on $A_\mathcal{U}$. In
this situation, we try to answer the following:

\begin{ques} \label{pregunta1}
When is it true that, for every $G$-graded algebra $A=\oplus_{g\in
G}A_g$, the subspace $A_\mathcal{U}=\oplus_{u\in\mathcal{U}}A_u$
inherites a structure of $G$-graded algebra with the
support-restricted multiplication?. For such a $\mathcal{U}$, how
are connected the graded module theories of $A$ and
$A_\mathcal{U}$?
\end{ques}

Also, concerning regradings, we try to answer the following:

\begin{ques} \label{pregunta2}
Let $G',G$ be groups and $\varphi :G'\longrightarrow G$ be a map.
When is it true that, for every $G$-graded algebra $B=\oplus_{g\in
G}B_g$ with support contained in $Im(\varphi )$, taking
$\tilde{B}_\sigma =B_{\varphi (\sigma )}$, for every $\sigma\in
G'$, gives rise to a $G'$-graded algebra $\tilde{B}$  which
coincides with $B$ as ungraded algebra?. In that case, how are
connected the graded module theories of $B$ and the (regraded)
$\tilde{B}$?
\end{ques}

  The structure of the paper goes as follows. In Section 2, we
  show that the
  answer to Question \ref{pregunta1} is: 'when $\mathcal{U}$ is a ring-supporting
  subset (see definition below) of $G$' (cf. Proposition
  \ref{caracterizacion-de-ring-supporting}).
We then introduce the concept of (right) modular pair of subsets
(see Definition \ref{definicion}) of $G$. If
$(\mathcal{S},\mathcal{U})$ is such a pair, in which case
$\mathcal{U}$ is ring-supporting, then the assigment
$M\rightsquigarrow M_\mathcal{S}$ gives an exact functor
$(-)_\mathcal{S}:Gr_A\longrightarrow Gr_{A_\mathcal{U}}$, which is
fully faithful when restricted to a well-described subcategory
$\mathcal{G}(\mathcal{S},\mathcal{U})$
 (Theorem \ref{fullyfaithful-killing-supports}). The essential
image of $\mathcal{G}(\mathcal{S},\mathcal{U})$ by
$(-)_\mathcal{S}$, denoted $\mathcal{L}(\mathcal{S},\mathcal{U})$,
contains all the  graded $A_\mathcal{U}$-modules presented in
degrees belonging to $(\mathcal{S}:\mathcal{U})=\{g\in G:$
$g\mathcal{U}=\mathcal{S}\}$. In Section 3, we show that the
answer to Question \ref{pregunta2} is: 'when $\varphi$ is a
pseudomorphism of groups'. Then we show that there is an
equivalence  between the category $\Sigma_{Im(\varphi )}$ of
$G$-graded $B$-modules $X$ with support $Supp(X)\subseteq
Im(\varphi )$ and the category $\tilde{\Sigma}_{Im(\varphi )}$ of
$G'$-graded $\tilde{B}$-modules $V$ such that
$V_\sigma\tilde{B}_\tau =0$ whenever $\varphi (\sigma )\varphi
(\tau )\not\in Im(\varphi )$ (cf. Proposition
\ref{caracterizacion-de-pseudomorphisms}). We then mix both
questions, namely, study the case when $\mathcal{U}$ is a
ring-supporting subset of $G$ which is the image of a
pseudomorphism $\varphi :G'\longrightarrow G$. Then, passing from
$A$ to $B=A_\mathcal{U}$ and from $B$ to $\tilde{B}$, we obtain an
equivalence of categories between the subcategory
$\mathcal{G}(\mathcal{S},\mathcal{U})$ of $Gr_A$ and a  full
subcategory $\mathcal{L}_{\tilde{B}}$ of $Gr_{\tilde{B}}$
(Proposition \ref{ring-supporting=image of pseudomorphism}).

In Section 4, we study the ring-supporting subsets
$\mathcal{U}\subseteq\mathbf{Z}$  of the integers, giving
essentially their structure (Theorem \ref{ring-supporting}). Then
precise criteria are given for a graded $A_\mathcal{U}$-module
 to
belong to $\mathcal{L}(\mathcal{S},\mathcal{U})$, criteria which
are specially simple when
 $\mathcal{U}$ is the translation of an inverval (cf. Corollary
 \ref{liftable-interval}). In addition, in this last case,
 $\mathcal{U}$ is the image of a precise pseudomorphism
 $\mathbf{Z}\longrightarrow\mathbf{Z}$. Then, for a
 positively  $\mathbf{Z}$-graded algebra $A=\oplus_{n\geq
 0}A_n$, taking $B=A_\mathcal{U}$ and then the regraded
 $\tilde{B}$, we obtain a precise description of the subcategories
 $\mathcal{G}(\mathcal{S},\mathcal{U})$ and
 $\mathcal{L}_{\tilde{B}}$ mentioned above. In case $\Lambda$ is a
 $n$-Koszul algebra and $A=\Lambda^!$, the results can be applied
 to show equivalences between appropriate full subcategories of
 $Gr_{\Lambda^!}$ and $Gr_E$ (cf. Remark \ref{n-Koszul}).

All algebras in the paper are associative with unit and modules
are right modules, unless explicitly said otherwise. Our basic
reference for graded rings and modules will be
 \cite{NVO}, while the terminology concerning torsion theories in Grothendieck categories can be found in \cite{St}.
 The few definitions and properties that we need about
  $n$-Koszul algebras can be found in \cite{B} and
 \cite{GMMVZ}. For the non-defined terms concerning algebras  and modules, we refer to
 \cite{ARS}.

\section{Killing of supports and induced equivalences}

In this and next section, we assume that $R$ is a commutative ring, $G$ is a group and $A=%
\underset{g\in G}{\oplus }A_{g}$ is a $G-$graded $R$-algebra.
 Let us fix a subset $\mathcal{S}$ of $G$
  and, for every graded $A$-module $M=\oplus_{g\in G}M_{g}$,  put
  $M_{\mathcal{S}}$ =$\oplus_{
s\in \mathcal{S}}M_{s}$,  letting $M_\emptyset=0$. We take
$\mathcal{T}=:\mathcal{T}_{\mathcal{S}}=\{ M\in Gr_A:$
$Supp(M)\subseteq G\setminus\mathcal{S} \}$, where $Supp(M)=\{g\in
G:$ $M_g\neq 0\}$ is the support of $M$. It is a hereditary
torsion class closed for products in $Gr_A$. We can view $R$ as a
$G$-graded ring, with $R_1=R$ and $R_g=0$ for $g\neq 1$, and then
an object of $Gr_R$ is just a family $(M_g)_{g\in G}$ of
$R$-modules. We view $M_\mathcal{S}$ as an object of $Gr_R$,  with
$(M_\mathcal{S})_g=0$ when $g\in G\setminus\mathcal{S}$.

\begin{prop} \label{faithful-killing-supports}
The assignment $M\longrightarrow M_{\mathcal{S}}$ gives an exact
functor, $(-)_{\mathcal{S}}:Gr_{A}\rightarrow Gr_R$ with kernel $
\mathcal{T}$.  Moreover, the induced functor
$Gr_{A}/\mathcal{T}\rightarrow Gr_R$ is faithful and identifies
$Gr_{A}/\mathcal{T}$ with a (not necessarily full) subcategory of
$Gr_R$.
\end{prop}

\begin{proof}
As a functor $Gr_{A}\longrightarrow Mod_R$ , $M\rightsquigarrow
M_{s}$,  is exact, for every $s\in \mathcal{S}$ . Consequently,
the assignment $M\rightsquigarrow M_{\mathcal{S}}$ yields an exact
functor: $Gr_{A}\rightarrow Gr_{R}$.

The kernel of $(-)_{\mathcal{S}}$ is $\{ M\in Gr_{A}:$
$M_{\mathcal{S}}=0\} =\mathcal{T} $.  Then the universal property
of the quotient category (cf \cite{G}) yields a unique $R$-linear
functor $F:Gr_{A}/\mathcal{T}\longrightarrow Gr_R$ such that the
diagram:

\vspace*{0.5cm}

\setlength{\unitlength}{1mm}
\begin{picture}(140,30)

\put(25,24){$Gr_A$} \put(32,25){\vector(1,0){53}}
\put(87,24){$Gr_R$} \put(55,2){$Gr_A/\mathcal{T}$}
\put(28,22){\vector(3,-2){26}} \put(62,5){\vector(3,2){25}}
\put(57,27){$(-)_\mathcal{S}$} \put(38,12){$p$} \put(79,13){$F$}
\end{picture}

is commutative, where $p$ is the canonical projection. All we need
to check is that $F$ is faithful. By \cite{G},  we have
$Hom_{Gr_{A}/ \mathcal{T}}(M,N)=\underset{N^{\prime
}\text{,}M/M'\in
\mathcal{T}}{\underrightarrow{\lim }}%
Hom_{Gr_{A}}(M^{\prime },N/N^{\prime }).$
 Here the order in the
index set is: $(M^{\prime },N^{\prime })\leq (M^{\prime \prime
},N^{\prime \prime })$ iff $M^{\prime \prime }\subseteq M^{\prime
}$ and $N^{\prime }\subseteq N^{\prime \prime }.$ In our situation
we have a unique maximum in that set. Indeed, since $\mathcal{T}$
is closed for products and subobjects in $Gr_A$, we have that
$\underset{M^{\prime }\leq M,M/M^{\prime }\in \mathcal{T}}{\cap
}M^{\prime }=\overset{\sim }{M}$ is a graded submodule of $M$ such
that $M/\overset{\sim }{M}\in \mathcal{T}.$ On the other hand,
$N^{\prime }\subseteq t(N),$ where $t(N)$ is by definition the
largest submodule of $N$ belonging to $\mathcal{T}.$

Therefore $(M^{\prime },N^{\prime })\leq (\overset{\sim
}{M},t(N))$ for all $(M^{\prime },N^{\prime })$ in the index
set.\newline
Then Hom$_{Gr_{A}/\mathcal{T}}($M$,$N$)$=Hom$%
_{Gr_{A}}$($\overset{\sim }{\text{M}}$,N/t(N))$.$

In order to prove the faithful condition of $F$ we need a more
precise
description of $\overset{\sim }{M}$ and $t(N).$ Notice that if \ $%
M^{\prime }\leq M$ is a graded submodule, then  $M/M^{\prime }\in
\mathcal{T}$ iff $M_{s}=M_{s}^{\prime }$ for every $s\in
\mathcal{S}$ iff  $M_{\mathcal{S}}\subseteq M^{\prime
}.$ As a consequence, $\overset{\sim }{M}=M_{\mathcal{S}}A$ is the graded $%
A- $ submodule of $M$ generated by $M_{\mathcal{S}}.$ On the other hand, in case $%
N^{\prime }\leq N$, we have that  $N^{\prime }\in \mathcal{T}$ iff
 $N_{s}^{\prime }=0$ for all $s\in \mathcal{S}.$ As a
consequence, $t(N)$ is the largest graded submodule of $N$
contained in $N_{G\smallsetminus \mathcal{S}}.$ In particular,
this implies that $N_{\mathcal{S}}=(N/$ $t(N)_{\mathcal{S}}$.

Now the canonical map
$Hom_{Gr_{A}}(M_{\mathcal{S}}A,N/t(N))\rightarrow
Hom_{Gr_R}(M_{\mathcal{S}},N_{\mathcal{S}})$
maps $f\in Hom_{Gr_{A}}(M_{\mathcal{S}}A,N/t(N))$ onto $%
f_{\mid M_{\mathcal{S}}}.$ If $F(f)=0,$ then $f_{\mid
M_{\mathcal{S}}}=0$ and, since $M_{\mathcal{S}}A$ is generated by
$M_{\mathcal{S}}$ as A- module, $f=0$. That proves the faithful
condition of $F$.

\end{proof}

An arbitrary killing of supports won't lead to a reasonably good
module theory over a new graded algebra. We need to impose some
restrictions.

\begin{defi} \label{definicion}
Let $G$ be a group and $(\mathcal{S},\mathcal{U})$ a pair of
nonempty subsets of $G$ such that $1\in\mathcal{U}$. We shall say
that $(\mathcal{S},\mathcal{U})$ is a  {\bf right premodular pair}
if whenever
$(s,u,v)\in\mathcal{S}\times\mathcal{U}\times\mathcal{U}$ is a
triple, with $suv\in\mathcal{S}$, one has that $su\in\mathcal{S}$
if, and only if, $uv\in\mathcal{U}$. By symmetry we get the notion
of left premodular pair.

A subset $\mathcal{U}\subseteq G$ containing $1$ will be called
{\bf ring-supporting} in case $(\mathcal{U},\mathcal{U})$ is a
(left or right) premodular pair. A  {\bf right} (resp. {\bf left})
{\bf modular pair} is a right (resp. left) premodular pair  such
that $\mathcal{U}$ is ring-supporting.
\end{defi}

The terms above are justified by the following results.

\begin{prop} \label{caracterizacion-de-ring-supporting}
Let $\mathcal{U}\subseteq G$ be a  subset of the group $G$. The
following assertions are equivalent:

\begin{enumerate}
\item $\mathcal{U}$ is ring-supporting \item For every $G$-graded
algebra $A=\oplus_{g\in G}A_g$, the sub-bimodule
$A_\mathcal{U}=\oplus_{u\in\mathcal{U}}A_u$ has a  structure of
$G$-graded $R$-algebra by putting $(A_\mathcal{U})_g=0$, when
$g\in G\setminus\mathcal{U}$, and considering the
support-restricted multiplication.

\end{enumerate}
\end{prop}

\begin{proof}
$1)\Longrightarrow 2)$ The ring-supporting condition of
$\mathcal{U}$ guarantees that the given multiplication in
$A_\mathcal{U}$ is associative and, from that, it follows easily
that $A_\mathcal{U}$ becomes a $G$-graded algebra.

$2)\Longrightarrow 1)$ Take the group algebra $A=RG$. Then
$A_\mathcal{U}=R\mathcal{U}=\oplus_{u\in\mathcal{U}}Ru$. Denote by
$\cdot$ the support-restricted multiplication on $R\mathcal{U}$.
Since $R\mathcal{U}$ has an identity, necessarily
$1\in\mathcal{U}$. On the other hand, in case
$u,v,w\in\mathcal{U}$ and $uvw\in\mathcal{U}$, one has $(u\cdot
v)\cdot w=uvw$ exactly when $uv\in\mathcal{U}$, and $uvw=u\cdot
(v\cdot w)$ exactly when $vw\in\mathcal{U}$. The associative
property of $\cdot$ then gives that $uv\in\mathcal{U}$ if, and
only if, $vw\in\mathcal{U}$.
\end{proof}

\begin{prop} \label{caracterizacion-de-right-modular}

Let $(\mathcal{S},\mathcal{U})$ be a pair of subsets of the group
$G$, where $\mathcal{U}$ is ring-supporting. The following
assertions are equivalent:

\begin{enumerate}
\item $(\mathcal{S},\mathcal{U})$ is right modular \item  For
every $G$-graded algebra and every $M\in Gr_A$, the sub-bimodule $M_{\mathcal{S}%
}$ =$\underset{s\in \mathcal{S}}{\oplus }M_{s}$ has a  structure
of graded right $A_{\mathcal{U}}-$ module by putting
$(M_{\mathcal{S}})_g=0$, when $g\in G\setminus\mathcal{S}$, and
defining the exterior multiplication of nonzero homogeneous
elements by: $x\cdot a=\left\{
\begin{array}{cc}
xa & \text{if }\deg (x)\text{deg( }a\text{)}\in \mathcal{S} \\
0 & \text{otherwise}%
\end{array}%
\right. .$

\end{enumerate}

In the above situation, the functor $(-)_\mathcal{S}$ of
Proposition \ref{faithful-killing-supports} may  be viewed as an
$R$-linear functor $Gr_A\longrightarrow Gr_{A_\mathcal{U}}$.

\end{prop}

\begin{proof}
  Along the lines of the previous proposition.
\end{proof}

We now give some properties of the just introduced subsets of a
group. In case $(\mathcal{S},\mathcal{U})$ is a  pair of subsets
of the group $G$, we shall put $(\mathcal{S}:\mathcal{U})=\{g\in
G:$ $g\mathcal{U}=\mathcal{S}\}$. If $1\in\mathcal{U}$ then
$(\mathcal{S}:\mathcal{U})$ is a (possibly empty) subset of
$\mathcal{S}$.

\begin{prop} \label{ejemplos-modular-pairs}
Let $G$ be a group and $\mathcal{S},\mathcal{U}$ be nonempty
subsets of $G$. The following assertions hold:
\begin{enumerate}
\item If $\mathcal{U}$ is  ring-supporting, then
$\mathcal{U}^{-1}$ is also  ring-supporting  \item The pair
$(\mathcal{S},\mathcal{U})$ is right modular if, and only if,
$(\mathcal{U}^{-1},\mathcal{S}^{-1})$ is left modular. \item
$(\mathcal{U},\mathcal{S})$ is  left modular  if, and only if,
$(\mathcal{S}^{-1},\mathcal{U})$ is  right modular  \item If
$\mathcal{U}$ is  ring-supporting  then
$(\mathcal{U},\mathcal{U}^{-1})$ is  left and right modular \item
If $1\in \mathcal{U}$  then $(\mathcal{U}:\mathcal{U})$ is a
subgroup of $G$ contained in $\mathcal{U}$ which, in case
$\mathcal{U}$ is ring-supporting, coincides with
$\mathcal{U}\cap\mathcal{U}^{-1}$ \item If
$(\mathcal{S},\mathcal{U})$ is a right modular pair and $g\in
(\mathcal{S}:\mathcal{U})$, then $\mathcal{S}=g\mathcal{U}$ and
$(\mathcal{S}:\mathcal{U})=g(\mathcal{U}:\mathcal{U})$
\end{enumerate}
\end{prop}

\begin{proof}
We leave to the reader the easy verification of 1 and 2. In 3, by
symmetry it will  be enough to check one of the implications.
Suppose then that $(\mathcal{U},\mathcal{S})$ is a left modular
pair.  We consider elements $s\in\mathcal{S}$ and
$u,v\in\mathcal{U}$ such that $s^{-1}uv\in\mathcal{S}^{-1}$ or,
equivalently, $v^{-1}u^{-1}s\in\mathcal{S}$. We want to prove that
$uv\in\mathcal{U}$ if, and only if, $s^{-1}u\in\mathcal{S}^{-1}$.
That is, we want to prove that $uv\in\mathcal{U}$ if, and only if,
$u^{-1}s\in\mathcal{S}$. But we have elements  $u$, $v$ in
$\mathcal{U}$ and $v^{-1}u^{-1}s\in\mathcal{S}$ with product
$uvv^{-1}u^{-1}s=s\in\mathcal{S}$. The left modularity of
$(\mathcal{U},\mathcal{S})$ gives that $uv\in\mathcal{U}$ if, and
only if, $vv^{-1}u^{-1}s=u^{-1}s\in\mathcal{S}$, as desired.

Assertion 4 follows from 2 and 3.

We prove now assertion 5. Clearly, $(\mathcal{U}:\mathcal{U})$ is
a subgroup of $G$ contained in $\mathcal{U}$. Then
$(\mathcal{U}:\mathcal{U})\subseteq\mathcal{U}\cap\mathcal{U}^{-1}$.
For the reverse inclusion in case $\mathcal{U}$ is
ring-supporting, take $u\in\mathcal{U}\cap \mathcal{U}^{-1}$ and
suppose that $u\not\in (\mathcal{U}:\mathcal{U})$. Then we have
two possibilities: i) there exists a $v\in\mathcal{U}$ such that
$uv\not\in\mathcal{U}$; ii) there exists $w\in\mathcal{U}$ such
that $w\not\in u\mathcal{U}$. In case i)  $u^{-1},u,v$ are
elements of $\mathcal{U}$ with product $u^{-1}uv=v\in\mathcal{U}$
such that $u^{-1}u=1\in\mathcal{U}$, but $uv\not\in\mathcal{U}$.
That contradicts the ring-supporting condition of $\mathcal{U}$.
In case ii), we have that $u^{-1}w\not\in\mathcal{U}$ and again we
get elements $u,u^{-1},w\in\mathcal{U}$ with product
$uu^{-1}w=w\in\mathcal{U}$ such that $uu^{-1}=1\in\mathcal{U}$,
but $u^{-1}w\not\in\mathcal{U}$. That is a contradiction. Then $u$
cannot exist, so that $\mathcal{U}\cap \mathcal{U}^{-1}\subseteq
(\mathcal{U}:\mathcal{U})$ as desired.

The proof of assertion 6 is left to the reader.
\end{proof}

\begin{lema} \label{torsionfree-cogeneration}
Let $A=\oplus_{g\in G}A_g$ be a $G$-graded algebra,
$\mathcal{S}\subseteq G$ be a subset and put
$\mathcal{T}=\mathcal{T}_\mathcal{S}$. The following assertions
are equivalent for a graded $A$-module $M$:

\begin{enumerate}
\item $M$ is $\mathcal{T}$-torsionfree (i.e. $Hom_{Gr_A}(T,M)=0$,
for all $T\in\mathcal{T}$) \item For every nonzero graded
submodule $N$ of $M$, one has $N_\mathcal{S}\neq 0$
\end{enumerate}
\end{lema}

\begin{proof}
If $T\in\mathcal{T}$ then $Supp(T)\subseteq
G\setminus\mathcal{S}$. Hence $M$ is $\mathcal{T}$-torsionfree if,
and only if, for every nonzero graded submodule $N$ of $M$, one
has $Supp(N)\not\subseteq G\setminus\mathcal{S}$.
\end{proof}

If $M\in Gr_A$ satisfies the conditions of the above lemma, we
shall indistinctly say that $M$ is $\mathcal{T}$-torsionfree or
that $M$ is {\bf cogenerated in degrees belonging to}
$\mathcal{S}$. In case $(\mathcal{S}, \mathcal{U})$ is a right
modular pair of subsets of $G$, the induced functor $F:Gr_{A}/
\mathcal{T}\longrightarrow Gr_{A_\mathcal{U}}$ is not full. We
shall see that it is full when restricted to an appropriate full
subcategory of $Gr_{A}/\mathcal{T}$.

\begin{lema} \label{posterior}
Let $(\mathcal{S},\mathcal{U})$ be a right modular pair of subsets
of the group $G$ and $A=\oplus_{\sigma\in G}A_\sigma$ be a
$G$-graded algebra. Suppose that $g\in (\mathcal{S}:\mathcal{U})$.
Then we have a commutative diagram of functors:

 \vspace*{0.5cm}

\setlength{\unitlength}{1mm}
\begin{picture}(140,30)

\put(42,25){$Gr_A$}\put(50,26){\vector(1,0){15}}
\put(68,25){$Gr_{A}$}

\put(42,4){$Gr_{A_\mathcal{U}}$}\put(50,5){\vector(1,0){15}}
\put(69,4){$Gr_{A_\mathcal{U}}$}

  \put(45,23){\vector(0,-1){14}}
\put(72,23){\vector(0,-1){14}}

\put(36,16){$(-)_{\mathcal{U}}$}  \put(73,16){$(-)_{\mathcal{S}}$}

\put(58,27){$\cong$}   \put(58,2){$\cong$}

\end{picture}

where the horizontal arrows are the shifting equivalences
$?[g^{-1}]:  V\rightsquigarrow V[g^{-1}]$
\end{lema}

\begin{proof}
We just need to prove that
$M[g^{-1}]_\mathcal{S}=M_\mathcal{U}[g^{-1}]$, for every $M\in
Gr_A$. Indeed $Supp(M_\mathcal{U}[g^{-1}])=\{h\in G:$ $g^{-1}h\in
Supp(M_\mathcal{U})\}$, which is clearly contained in
$g\mathcal{U}=\mathcal{S}$. Now if $s\in\mathcal{S}$ then
$M[g^{-1}]_s=M_{g^{-1}s}$ while
$M_\mathcal{U}[g^{-1}]_s=(M_\mathcal{U})_{g^{-1}s}$. But
$(M_\mathcal{U})_{g^{-1}s}=M_{g^{-1}s}$ since
$g^{-1}s\in\mathcal{U}$.

\end{proof}

 In order to give our next result we
consider the full subcategory
$\mathcal{G}(\mathcal{S},\mathcal{U})$ of $Gr_A$ with objects the
graded $A$-modules generated in degrees belonging to
$(\mathcal{S}:\mathcal{U})$ and cogenerated in degrees belonging
to $\mathcal{S}$ (i.e. $\mathcal{T}_\mathcal{S}$-torsionfree). On
the other hand, the objects of the essential image of the functor
$(-)_{\mathcal{S}}:Gr_A\longrightarrow Gr_{A_{\mathcal{U}}}$ will
be called {\bf liftable} graded $A_{\mathcal{U}}$-modules with
respect to $(-)_\mathcal{S}$. We shall denote by
$\mathcal{L}(\mathcal{S},\mathcal{U})$ the full subcategory of
$Gr_{A_\mathcal{U}}$ with objects the liftable
$A_\mathcal{U}$-modules generated in degrees belonging to
$(\mathcal{S}:\mathcal{U})$. A graded module (over any $G$-graded
algebra) $M$ will be called {\bf presented in degrees belonging
to} $\mathcal{X}\subseteq G$, when there is an exact sequence of
graded modules $P\longrightarrow Q\twoheadrightarrow M\rightarrow
0$, with $P,Q$ gr-projective and generated in degrees belonging to
$\mathcal{X}$.  We are now ready to prove:

\begin{teor} \label{fullyfaithful-killing-supports}
Let $(\mathcal{S},\mathcal{U})$ be a right modular pair of subsets
of $G$ and $\mathcal{T}=\mathcal{T}_\mathcal{S}$. The following
assertions hold:

\begin{enumerate}
\item The composition
$\mathcal{G}(\mathcal{S},\mathcal{U})\hookrightarrow
Gr_A\stackrel{p}{\longrightarrow}Gr_A/\mathcal{T}$ is fully
faithful \item The restriction of $(-)_\mathcal{S}$ to
$\mathcal{G}(\mathcal{S},\mathcal{U})$ is a fully faithful functor
which induces an equivalence of categories
$\mathcal{G}(\mathcal{S},\mathcal{U})\stackrel{\cong}{\longrightarrow}
\mathcal{L}(\mathcal{S},\mathcal{U})$ \item For every $g\in
(\mathcal{S}:\mathcal{U})$, there is a commutative diagram of
equivalences of categories:

 \vspace*{0.5cm}

\setlength{\unitlength}{1mm}
\begin{picture}(140,30)

\put(38,25){$\mathcal{G}(\mathcal{U},\mathcal{U})$}\put(50,26){\vector(1,0){15}}
\put(68,25){$\mathcal{G}(\mathcal{S},\mathcal{U})$}

\put(38,4){$\mathcal{L}(\mathcal{U},\mathcal{U})$}\put(50,5){\vector(1,0){15}}
\put(69,4){$\mathcal{L}(\mathcal{S},\mathcal{U})$}

  \put(45,23){\vector(0,-1){14}}
\put(72,23){\vector(0,-1){14}}

\put(36,16){$(-)_{\mathcal{U}}$} \put(46,16){$\cong$}
  \put(73,16){$(-)_{\mathcal{S}}$}

\put(58,27){$\cong$} \put(69,16){$\cong$}  \put(58,2){$\cong$}

\end{picture}

where the horizontal arrows are the shifting equivalences
$?[g^{-1}]$
\end{enumerate}

Moreover, $\mathcal{L}(\mathcal{S},\mathcal{U})$ contains all the
graded $A_\mathcal{U}$-modules presented in degrees belonging to
$(\mathcal{S}:\mathcal{U})$
\end{teor}

\begin{proof} Let us take $M,N\in\mathcal{G}(\mathcal{S},\mathcal{U})$. By the proof
of Proposition \ref{faithful-killing-supports},  we know that
$Hom_{Gr_A/\mathcal{T}}(M,N)=Hom_{Gr_A}(\tilde{M},N/t(N))$, where
$\tilde{M}$ is the $A$-submodule of $M$ generated by
$M_\mathcal{S}$. Since
$(\mathcal{S}:\mathcal{U})\subseteq\mathcal{S}$ and $M$ is
generated in degrees belonging to $(\mathcal{S}:\mathcal{U})$,  we
conclude that $\tilde{M}=M$. On the other hand, since $N$ is
torsionfree we get $t(N)=0$ and, hence,
$Hom_{Gr_A/\mathcal{T}}(M,N)=Hom_{Gr_A}(M,N)$, which proves
assertion 1.

By assertion 1, we can view $\mathcal{G}(\mathcal{S},\mathcal{U})$
as a full subcategory of $Gr_A/\mathcal{T}$, and then  the
restriction of $(-)_\mathcal{S}$ to
$\mathcal{G}(\mathcal{S},\mathcal{U})$, which we denote by $H$ in
the sequel,  can be identified with the restriction of
$F:Gr_A/\mathcal{T}\longrightarrow Gr_{A_\mathcal{U}}$ to
$\mathcal{G}(\mathcal{S},\mathcal{U})$. Since $F$ is faithful (cf.
Proposition \ref{faithful-killing-supports}) $H$ is faithful. We
next prove that it is full. Let
$M,N\in\mathcal{G}(\mathcal{S},\mathcal{U})$ and
$M_{\mathcal{S}}\overset{g}{\rightarrow }N_{\mathcal{S}}$ be a
morphism in $Gr_{A_{\mathcal{U}}}.$
Notice that  we have epimorphisms $\underset{i\in I}{%
\oplus }A\left[ s_{i}^{-1}\right] $ $\overset{\varepsilon _{M}}{\rightarrow }%
M$ and $\underset{j\in J}{\oplus }A\left[ t_{j}^{-1}\right] \overset{%
\varepsilon _{N}}{\rightarrow }N$ in $Gr_{A},$ where
$(s_{i})_{i\in I}$ and $(t_{j})_{j\in J}$ are families in
$(\mathcal{S}:\mathcal{U}).$ From the proof of Lemma
\ref{posterior} we get that
$A[g^{-1}]_\mathcal{S}=A_\mathcal{U}[g^{-1}]$, for every $g\in
(\mathcal{S}:\mathcal{U})$. Then, by applying $(-)_{\mathcal{S}}$,  we get epimorphisms $\underset{%
i\in I}{\oplus }A_{\mathcal{U}}\left[ s_{i}^{-1}\right] $ $\overset{\overset{%
\sim }{\varepsilon }_{M}}{\rightarrow }M_{\mathcal{S}}$ and $\underset{j\in J%
}{\oplus }A_{\mathcal{U}}\left[ t_{j}^{-1}\right] \overset{\overset{\sim }{%
\varepsilon }_{N}}{\rightarrow }N_{\mathcal{S}}$ in
$Gr_{A_\mathcal{U}}$
 By the projective condition of $\underset{i\in I}{\oplus }A_{\mathcal{U}}\left[ s_{i}^{-1}%
\right] $ in $Gr_{A_\mathcal{U}}$ we get a commutative diagram in
this latter category:

\begin{center}
$%
\begin{array}{ccc}
\underset{i\in I}{\oplus }A_{\mathcal{U}}\left[ s_{i}^{-1}\right] & \overset{%
\overset{\sim }{\varepsilon }_{M}}{\rightarrow } & M_{\mathcal{S}} \\
\downarrow \phi &  & \downarrow g \\
\underset{j\in J}{\oplus }A_{\mathcal{U}}\left[ t_{j}^{-1}\right] & \overset{%
\overset{\sim }{\varepsilon }_{N}}{\rightarrow } & N_{\mathcal{S}}%
\end{array}%
$
\end{center}

Recall that $\underset{i\in I}{\oplus }A_{\mathcal{U}}\left[ s_{i}^{-1}%
\right] =A_{\mathcal{U}}^{(I)}$ as ungraded right
$A_{\mathcal{U}}-$ module,
with the grading obtained by assigning degree $s_{i}$ to the i-th element $%
e_{i}$ of the canonical basis of $A_{\mathcal{U}}^{(I)}$ as
$A_{\mathcal{U}}
$ - module. Then $\phi (e_{i})$ is an element of $\underset{j\in J}{\oplus }%
A_{\mathcal{U}}\left[ t_{j}^{-1}\right] $ of degree $s_{i}.$ i.e.
$\phi
(e_{i})$ $\in $ $\underset{j\in J}{\oplus }A_{\mathcal{U}}\left[ t_{j}^{-1}%
\right]_{s_i} =\underset{j\in J}{\oplus }A_{t_{j}^{-1}s_{i}}.$ We define now $%
\overset{\_}{\phi }:\underset{i\in I}{\oplus }A\left[ s_{i}^{-1}\right] $ $%
\rightarrow $ $\underset{j\in J}{\oplus }A\left[ t_{j}^{-1}\right]
$ as the
only homomorphism of (graded) right $A$-modules such that $\overset{\_}{\phi }%
(e_{i})=\phi (e_{i})$ for all $i.$

This yields a diagram in $Gr_{A}:$

\begin{center}
$%
\begin{array}{ccc}
\underset{i\in I}{\oplus }A\left[ s_{i}^{-1}\right] &
\overset{\varepsilon
_{M}}{\rightarrow } & M \\
\downarrow \overset{\_}{\phi } &  &  \\
\underset{j\in J}{\oplus }A\left[ t_{j}^{-1}\right] & \overset{\varepsilon _{N}}%
{\rightarrow } & N%
\end{array}%
$
\end{center}

We claim that  $\varepsilon _{N}\overset{\_}{\phi }$ vanishes on $%
Ker\varepsilon _{M},$ which implies the existence of a unique morphism $%
f:M\rightarrow N$ in $Gr_{A}$ such that $\varepsilon _{N}%
\overset{\_}{\phi }$ = $f\varepsilon _{M}.$ From that one easily
gets that $f_\mathcal{S}=g$ and the full condition for $H$ will
follow. In order to prove our claim, let $x\in Ker\epsilon _{M}$
be a homogeneous element. If $\deg (x)\in \mathcal{S}$, then $x\in
Ker(\overset{\sim}{\varepsilon}_M)$, so that
$\varepsilon_N\bar{\phi}(x)=\overset{\sim}{\varepsilon}_N\phi
(x)=g\overset{\sim}{\varepsilon}_M(x)=0$. In case $\deg (x)\notin
\mathcal{S}$,
we have $\deg ((\varepsilon _{N}\overset{\_}{\phi }%
)(x))=\deg (x)\notin \mathcal{S}$. In either case,  every element
in
$(\varepsilon _{N}\overset{%
\_}{\phi })(Ker\varepsilon _{M})$ has support in $G\smallsetminus
\mathcal{S} $ . This implies $(\varepsilon _{N}\overset{\_}{\phi
})(Ker\varepsilon _{M})\subseteq t(N)=0$  as desired.

In order to end the proof of assertion 2, it only remains to prove
that if $X\in\mathcal{L}(\mathcal{S},\mathcal{U})$  then there is
a $M\in Gr_A$ such that $M_\mathcal{S}=X$ and $M$ is generated in
degrees belonging to $(\mathcal{S}:\mathcal{U})$. In that case, we
could always assume that $M$ is $\mathcal{T}$-torsionfree, so that
$M\in\mathcal{G}(\mathcal{S},\mathcal{U})$. Since $X$ is liftable,
we can  fix a $N\in Gr_A$ such that $N_\mathcal{S}=X$. We take a
family $(x_i)_{i\in I}$ of homogeneous generators of $X$, where
$deg(x_i)=:s_i\in (\mathcal{S}:\mathcal{U})$ for all $i\in I$.
Then $x_i\in X_{s_i}=N_{s_i}$ and we can consider $M=\sum_{i\in
I}x_iA$, i.e., the graded $A$-submodule of $N$ generated by the
$x_i$. Then $X=\sum_{i\in I}x_iA_\mathcal{U}\subseteq
M_\mathcal{S}$, due to the fact that $s_i\in
(\mathcal{S}:\mathcal{U})$ for all $i\in I$. From that we get that
$X\subseteq M_\mathcal{S}\subseteq N_\mathcal{S}=X$ and, hence,
equalities. Then $M_\mathcal{S}=X$ as desired.

Assertion 3 follows directly from Lemma \ref{posterior}. For the
final statement, notice that every gr-projective
$A_\mathcal{U}$-module $P$ generated in degrees belonging to
$(\mathcal{S}:\mathcal{U})$ is a direct summand of $\oplus_{i\in
I}A_\mathcal{U}[g_i^{-1}]$, where $(g_i)_{i\in I}$ is a family of
elements of $(\mathcal{S}:\mathcal{U})$. The proof of Lemma
\ref{posterior} implies that $P$ is liftable and, hence, belongs
to $\mathcal{L}(\mathcal{S},\mathcal{U})$. But, due to assertion 2
and the exactness of $(-)_\mathcal{S}$, the class
$\mathcal{L}(\mathcal{S},\mathcal{U})$ is closed for cokernels.
\end{proof}

\begin{ejem}
Let $G$ be a group, $A=\oplus_{g\in G}A_g$ be a $G$-graded algebra
and $H<G$ be a subgroup. We put $\mathcal{S}=\mathcal{U}=H$. Then
$A_H=\oplus_{h\in H}A_h$ is in this case a subalgebra of $A$.
Moreover, $(\mathcal{S}:\mathcal{U})=(H:H)=\{g\in G:$ $gH=H\}=H$.
In this case, $\mathcal{G}(H,H)$ is the full subcategory of $Gr_A$
with objects those graded $A$-modules $M$ such that $M=M_HA$ and
$M$ contains no nonzero graded submodule with support in
$G\setminus H$.  Our Theorem \ref{fullyfaithful-killing-supports}
says that the functor $(-)_H:Gr_A\longrightarrow Gr_{A_H}$
establishes equivalences of categories
$\mathcal{G}(H,H)\stackrel{\cong}{\longrightarrow}
Gr_A/\mathcal{T}\stackrel{\cong}{\longrightarrow}Gr_{A_H}$,
understanding $Gr_{A_H}$ as the category of $G$-graded
$A_H$-modules with support in $H$. The reason is that every graded
$A_H$-module $X$ is liftable, because $M=:X\otimes_{A_H}A$ is a
graded $A$-module generated in degrees belonging to $H$ such that
$M_H\cong X$.
\end{ejem}

\section{Pseudomorphisms of groups and regradings}

In this section we answer Question \ref{pregunta2} (see the
introduction).

\begin{defi}
Let $G,G'$ be groups. A {\bf pseudomorphism} $\varphi
:G'\longrightarrow G$ is an injective map satisfying the following
two properties:

\begin{enumerate}
\item $\varphi (1)=1$ \item If $\sigma ,\tau\in G'$ and
 $\varphi (\sigma )\varphi (\tau )\in Im(\varphi )$, then
 $\varphi (\sigma\tau )=\varphi (\sigma )\varphi (\tau )$
 \end{enumerate}
\end{defi}

\begin{prop} \label{caracterizacion-de-pseudomorphisms}
Let $G',G$ be groups and $\varphi :G'\longrightarrow G$ be a map. The following
assertions are equivalent:

\begin{enumerate}
\item[a)] $\varphi$ is a pseudomorphism of groups \item[b)]  $1\in
Im(\varphi )$ and, for every $G$-graded algebra $B=\oplus_{g\in
G}B_g$ with support contained in $Im(\varphi )$, putting
$\tilde{B}=B$ and $\tilde{B}_\sigma =B_{\varphi (\sigma )}$, for
every $\sigma\in G'$, one obtains a $G'$-graded algebra
$\tilde{B}$ coinciding with $B$ as an ungraded algebra
\end{enumerate}

If $B$ and $\tilde{B}$ are $G$-graded and $G'$-graded algebras as
in b), then, for every element $g\in G$, there is an equivalence
of (abelian) categories between:

\begin{enumerate}
\item The full subcategory $\Sigma_{gIm(\varphi )}$ of $Gr_B$ with
objects the graded $B$-modules $X$ such that $Supp(X)\subseteq
gIm(\varphi )$ \item  The full subcategory
$\tilde{\Sigma}_{Im(\varphi )}$ of $Gr_{\tilde{B}}$ with objects
the graded $\tilde{B}$-modules $V$ such that
$V_\sigma\tilde{B}_\tau=0$ whenever $\varphi (\sigma )\varphi
(\tau )\not\in Im(\varphi )$
\end{enumerate}
\end{prop}

\begin{proof}
$a)\Longrightarrow b)$ Let $B=\oplus_{g\in G}B_g$ be any
$G$-graded algebra with $Supp(B)\subseteq Im(\varphi )$, and put
$\tilde{B}=B$ and $\tilde{B}_\sigma =B_{\varphi (\sigma )}$, for
all $\sigma\in G'$. Then
$\tilde{B}_\sigma\tilde{B}_\tau=B_{\varphi (\sigma )}B_{\varphi
(\tau )}$ is  zero,  in case $\varphi (\sigma )\varphi (\tau
)\not\in Im(\varphi )$, and is an $R$-submodule of $B_{\varphi
(\sigma\tau )}=\tilde{B}_{\sigma\tau}$, in case $\varphi (\sigma
)\varphi (\tau )\in Im(\varphi )$. Hence we always have
$\tilde{B}_\sigma\tilde{B}_\tau\subseteq\tilde{B}_{\sigma\tau}$,
so that $\tilde{B}$ is a $G'$-graded algebra.

$b)\Longrightarrow a)$ Since $1\in Im(\varphi )$, by taking
$G$-graded algebras concentrated in degree $1$, we can ensure that
there exist graded $G$-algebras with support contained in
$Im(\varphi )$. Moreover, for any such algebra $B=\oplus_{g\in
G}B_g$, the hypothesis says that putting $\tilde{B}=B$ and
$\tilde{B}_\sigma =B_{\varphi (\sigma )}$, for every $\sigma\in
G'$, we get a $G'$-graded algebra $\tilde{B}$ coinciding with $B$
as an ungraded algebra. Since the identity of $\tilde{B}=B$ is
always a homogeneous element of degree $1$, we conclude that
$1\in\tilde{B}_1=B_{\varphi (1)}$ and $1\in B_1$. But, by
definition of $G$-grading, the sum $\oplus_{g\in G}B_g$ is direct.
It follows that $\varphi (1)=1$.

On the other hand, if $1\neq g\in Im(\varphi )$
  then we can always find a $G$-graded algebra $B$, with $Supp(B)\subseteq Im(\varphi
)$, such that $B_g\neq 0$. Indeed, take the algebra $B=R[x]/(x^2)$
with $G$-grading given by putting $B_1=R$ and $deg(\bar{x})=g$.
Suppose now that $\sigma ,\tau\in G'$ are elements such that
$\varphi (\sigma )=\varphi (\tau)$. Then we consider any
$G$-graded algebra $B$, with $Supp(B)\subseteq Im(\varphi )$, such
that $B_{\varphi (\sigma )}\neq 0$. Then we have that
$\tilde{B}_\sigma =B_{\varphi (\sigma )}=B_{\varphi (\tau
)}=\tilde{B}_\tau$ is nonzero.  Then, by definition of
$G'$-grading,  we necessarily have $\sigma =\tau$. That proves
that $\varphi$ is injective.

Finally,  if $\sigma ,\tau\in G'$ and $\varphi (\sigma )\varphi
(\tau )\in Im(\varphi )$, then we have
$\tilde{B}_\sigma\tilde{B}_\tau=B_{\varphi (\sigma )}B_{\varphi
(\tau )}\subseteq B_{\varphi (\sigma )\varphi (\tau )}$ and, since
$\tilde{B}$ is a $G'$-graded algebra, we also have
$\tilde{B}_\sigma\tilde{B}_\tau\subseteq\tilde{B}_{\sigma\tau}=B_{\varphi
(\sigma\tau )}$. It follows that $\varphi (\sigma\tau )=\varphi
(\sigma )\varphi (\tau )$ unless $B_{\varphi (\sigma )}B_{\varphi
(\tau )}=0$, for every graded $G$-algebra $B$ with
$Supp(B)\subseteq Im(\varphi )$. We should then discard this last
possibility, for which we prove that if $g,h\in G\setminus{1}$ are
elements such that $g,h,gh\in Im(\varphi )$, then there exists a
$G$-graded algebra $B$, with $Supp(B)\subseteq Im(\varphi )$, such
that $B_gB_h\neq 0$. We consider all the possibilities:

\begin{enumerate}
\item[i)] $g=h$ and $gh=1$:  Take the group algebra $B=RC$, where
$C=<g>$ is the subgroup of $G$ generated by $g$ with the obvious
$G$-grading  \item[ii)] $g=h$ and $gh\neq 1$: Take $B=R[x]/(x^3)$,
$B_1=R$ and $deg(\bar{x})=g$  \item[iii)] $g\neq h$ and $gh=1$:
Take $B=R<x,y>/(x^2,y^2,yx)$, with $B_1=R+R\bar{x}\bar{y}$,
$deg(\bar{x})=g$ and $deg(\bar{y})=h$ \item[iv)] $g\neq h$ and
$gh\neq 1$: Take $B=R<x,y>/(x^2,y^2,yx)$, with $B_1=R$,
$deg(\bar{x})=g$ and $deg(\bar{y})=h$
\end{enumerate}

For the final part, given a subset $\mathcal{X}\subseteq G$, we
denote by $\Sigma_\mathcal{X}$ the full subcategory of $Gr_B$ with
objects those $M\in Gr_B$ such that $Supp(M)\subseteq\mathcal{X}$.
It is clear that, for every $g\in G$, the shifting equivalence
$?[g]:Gr_B\stackrel{\cong}{\longrightarrow}Gr_B$ induces by
restriction an equivalence $\Sigma_{gIm(\varphi
)}\cong\Sigma_{Im(\varphi )}$. So it is not restrictive to assume,
something that we do in the rest of the proof, that $g=1$. Let
then  $X=\oplus_{g\in G}X_g$ be a graded $B$-module with
$Supp(X)\subseteq Im(\varphi )$. Then we define $\tilde{X}=X$ as
ungraded $\tilde{B}$- (or $B$-)module and give it a structure of
$G'$-graded $\tilde{B}$-module by putting $\tilde{X}_\sigma
=X_{\varphi (\sigma )}$, for all $\sigma\in G'$. It is a mere
routine to check that the assignment $X\rightsquigarrow\tilde{X}$
defines a fully faithful exact functor $F:\Sigma_{Im(\varphi
)}\longrightarrow Gr_{\tilde{B}}$. Moreover, if $\sigma ,\tau\in
G'$ are elements such that $\varphi (\sigma )\varphi (\tau
)\not\in Im(\varphi )$, in which case $X_{\varphi (\sigma )\varphi
(\tau )}=0$, then we have $\tilde{X}_\sigma\tilde{B}_\tau
=X_{\varphi (\sigma )}B_{\varphi (\tau )}\subseteq X_{\varphi
(\sigma )\varphi (\tau )}=0$, thus proving that
$Im(F)\subseteq\tilde{\Sigma}_{Im(\varphi )}$. Conversely, if
$V\in\tilde{\Sigma}_{Im(\varphi )}$ then we can form a $G$-graded
$B$-module $X$ as follows. We put $X=V$ as ungraded $B-$ (or
$\tilde{B}-$)module and define $X_g=0$, when $g\not\in Im(\varphi
)$, and $X_g=V_\sigma$ provided $g=\varphi (\sigma )$. The
injective condition of $\varphi$ implies that the given
$G$-grading on $X$ is well-defined. We leave to the reader the
routinary task of checking that $X\in\Sigma_{Im(\varphi )}$ and
$\tilde{X}=V$. From that it follows that $F$ induces the desired
equivalence $\Sigma_{Im(\varphi
)}\stackrel{\cong}{\longrightarrow}\tilde{\Sigma}_{Im(\varphi )}$.

\end{proof}

It turns out that some of the most interesting ring-supporting
subsets of $\mathbf{Z}$ are the image of a pseudomorphism (see
next section). It is then of interest to consider the case when
$\mathcal{U}\subseteq G$ is a ring-supporting subset such that,
for some pseudomorphism $\varphi :G'\longrightarrow G$, one has
$Im(\varphi )=\mathcal{U}$. In that situation, starting with any
$G$-graded algebra $A$, we can take $B=A_\mathcal{U}$ in the above
proposition, so that we get a $G'$-graded algebra
$\tilde{B}=\tilde{A}_\mathcal{U}$. If now
$(\mathcal{S},\mathcal{U})$ is a right modular pair of subsets of
$G$ such that $(\mathcal{S}:\mathcal{U})\neq\emptyset$ then,
according to Proposition \ref{ejemplos-modular-pairs}(6), we have
$\mathcal{S}=g\mathcal{U}=gIm(\varphi )$, where $g\in
(\mathcal{S}:\mathcal{U})$. Moreover, since the objects in the
image of the functor $(-)_\mathcal{S}:Gr_A\longrightarrow
Gr_{A_\mathcal{U}}$ are graded $A_\mathcal{U}$-modules with
support in $\mathcal{S}$, we get a well-defined $R$-lineal exact
functor $\Phi :Gr_A\longrightarrow Gr_{\tilde{B}}$, which is the
following composition
$Gr_A\stackrel{(-)_\mathcal{S}}{\longrightarrow}\Sigma_{\mathcal{S}}
\stackrel{\cong}{\longrightarrow}\tilde{\Sigma}_{\mathcal{U}}\hookrightarrow
Gr_{\tilde{B}}$, where the central equivalence is that of
Proposition \ref{caracterizacion-de-pseudomorphisms}. We now have:

\begin{prop} \label{ring-supporting=image of pseudomorphism}
In the above situation, the following assertions hold:

\begin{enumerate}
\item $H'=\varphi^{-1}[(\mathcal{U}:\mathcal{U})]$ is a subgroup
of $G'$ \item  The functor $\Phi$ induces an equivalence of
categories between $\mathcal{G}(\mathcal{S},\mathcal{U})$ and the
full subcategory $\mathcal{L}_{\tilde{B}}$ of $Gr_{\tilde{B}}$
consisting of those $V\in Gr_{\tilde{B}}$ such that
$V\cong\tilde{X}$, for some
$X\in\mathcal{L}(\mathcal{S},\mathcal{U})$.  \item All
$V\in\mathcal{L}_{\tilde{B}}$ are graded $\tilde{B}$-modules
generated in degrees belonging to $H'$ and
$\mathcal{L}_{\tilde{B}}$ contains all the graded
$\tilde{B}$-modules presented in degrees belonging to $H'$

\end{enumerate}
\end{prop}

\begin{proof}
1) We prove, more generally, that if $\varphi :G'\longrightarrow
G$ is a pseudomorphism and $H<G$ is a subgroup contained in
$Im(\varphi )$, then $H'=\varphi^{-1}(H)$ is a subgroup of $G'$.

If $\sigma ,\tau\in H'$ then $\varphi (\sigma ), \varphi (\tau
)\in H$. Since $H$ is a subgroup of $G$, we have that $\varphi
(\sigma )\varphi (\tau )\in H\subseteq Im(\varphi )$. The fact
that $\varphi$ is a pseudomorphism implies that $\varphi
(\sigma\tau )=\varphi (\sigma )\varphi (\tau )$ and, hence,
$\sigma\tau\in\varphi^{-1}(H)=H'$.

On the other hand, if $\sigma\in H'$ then $\varphi (\sigma
)^{-1}\in H\subseteq Im(\varphi )$ due to the fact that $H$ is a
subgroup of $G$. Then we can choose $\tau\in G'$ such that
$\varphi (\sigma )^{-1}=\varphi (\tau )$. But since $\varphi
(\sigma )\varphi (\tau )=1\in Im(\varphi )$, we get that $\varphi
(1)=1=\varphi (\sigma )\varphi (\tau )=\varphi (\sigma\tau )$. The
fact that $\varphi$ is injective implies that $1=\sigma\tau$. Then
$\tau =\sigma^{-1}$ and hence $\sigma^{-1}\in H'=\varphi^{-1}(H)$.
That proves that $H'$ is a subgroup of $G'$.

2) According to the definition of $\Phi$ and Theorem 2.7, we get
equivalences
$\mathcal{G}(\mathcal{S},\mathcal{U})\stackrel{\cong}{\longrightarrow}
\mathcal{L}(\mathcal{S},\mathcal{U})\stackrel{\cong}{\longrightarrow}\mathcal{L}_{\tilde{B}}$.

3) Put $B=A_\mathcal{U}$. Then, as shown in the proof of
Proposition \ref{caracterizacion-de-pseudomorphisms},  the
equivalence $\Sigma_\mathcal{S}=\Sigma_{gIm(\varphi
)}\stackrel{\cong}{\longrightarrow}\tilde{\Sigma}_{Im(\varphi)}=\tilde{\Sigma}_\mathcal{U}$
is a composition
$\Sigma_{\mathcal{S}}\stackrel{[g]}{\longrightarrow}\Sigma_{\mathcal{U}}
\stackrel{\cong}{\longrightarrow}\tilde{\Sigma}_\mathcal{U}$,
where the second equivalence takes graded $B$-modules generated
(resp. presented) in degrees belonging to
$(\mathcal{U}:\mathcal{U})$ to graded $\tilde{B}$-modules
generated (resp. presented) in degrees belonging to
$H'=\varphi^{-1}[(\mathcal{U}:\mathcal{U})]$. On the other hand,
the first equivalence
$[g]:\Sigma_{\mathcal{S}}\stackrel{\cong}{\longrightarrow}\Sigma_{\mathcal{U}}$
takes graded $B$-modules generated (resp. presented) in degrees
belonging to $(\mathcal{S}:\mathcal{U})$ to graded $B$-modules
generated (resp. presented) in degrees belonging to
$(\mathcal{U}:\mathcal{U})$. Now, from the chain of inclusions
 $\{X\in Gr_B:$ $X\text{ presented in
degrees belonging to }(\mathcal{S}:\mathcal{U})\}
\subseteq\mathcal{L}(\mathcal{S}:\mathcal{U})\subseteq\{X\in
Gr_B:$ $X\text{ generated in degrees belonging to
}(\mathcal{S}:\mathcal{U})\}\subseteq\Sigma_\mathcal{S}$ (cf.
Theorem \ref{fullyfaithful-killing-supports}),  assertion 3
follows.

\end{proof}

\section{The case $G=\mathbf{Z}$}

In this section we particularize the previous arguments to
$G=\mathbf{Z}$, the additive group of integers.  We thereby shift
from multiplicative to additive notation, and our algebra will be
a $\mathbf{Z}$-graded algebra $A=\oplus_{i\in\mathbf{Z}}A_i$. Our
first goal is to get as much information as possible  on
ring-supporting subsets $\mathcal{U}\subseteq\mathbf{Z}$ and
hence, by Proposition \ref{ejemplos-modular-pairs},    on modular
pairs $(\mathcal{S},\mathcal{U})$ such that $\emptyset\neq
(\mathcal{S}:\mathcal{U})=:\{m\in\mathbf{Z}:$
$m+\mathcal{U}=\mathcal{S}\}$. Those modular pairs will be of the
form $(m+\mathcal{U},\mathcal{U})$, for some integer $m$.

\begin{lema} \label{reduction-modulo-(U:U)}
Let $\mathcal{U}\subseteq\mathbf{Z}$ be a subset containing $0$
such that $(\mathcal{U}:\mathcal{U})=n\mathbf{Z}$, and let
$p:\mathbf{Z}\longrightarrow\mathbf{Z}_n$ be the canonical
projection. The following assertions are equivalent:

\begin{enumerate}
\item  $\mathcal{U}$ is a ring-supporting subset of $\mathbf{Z}$
\item  $\bar{\mathcal{U}}=p(\mathcal{U})$ is a ring-supporting
subset of $\mathbf{Z}_n$ such that
$(\bar{\mathcal{U}}:\bar{\mathcal{U}})=\{\bar{0}\}$
\end{enumerate}

In that case
$\mathcal{U}=\bigcup_{i\in\bar{\mathcal{U}}}(n\mathbf{Z}+i)$,
viewing $\bar{\mathcal{U}}$ as a subset of $[0,n)$ in case $n>0$.
\end{lema}

\begin{proof}

 We put $\bar{m}=p(m)$, for every $m\in\mathbf{Z}$. Let us first notice
 that  $p^{-1}(\bar{m}+\bar{\mathcal{U}}))=m+\mathcal{U}$. Indeed,
if $v\in p^{-1}(\bar{m}+\bar{\mathcal{U}}))$ then
$\bar{v}=p(v)=\bar{m}+\bar{u}$, for some $u\in\mathcal{U}$, which
implies that $v-m-u\in n\mathbf{Z}=(\mathcal{U}:\mathcal{U})$.
Then $v=m+u+(v-m-u)\in m+(v-m-u)+\mathcal{U}$. But
$(v-m-u)+\mathcal{U}=\mathcal{U}$, and then $v\in m+\mathcal{U}$
as desired. In particular $m+\mathcal{U}=\mathcal{U}$ iff
$\bar{m}+\bar{\mathcal{U}}=\bar{\mathcal{U}}$. In other words,
$m\in (\mathcal{U}:\mathcal{U})$ iff $\bar{m}\in
(\bar{\mathcal{U}}:\bar{\mathcal{U}})$. That gives that
$(\bar{\mathcal{U}}:\bar{\mathcal{U}})=\{\bar{0}\}$.

$1)\Longrightarrow 2)$ If  $u,v,w\in\mathcal{U}$ are elements such
that $\bar{u}+\bar{v}+\bar{w}\in\bar{\mathcal{U}}$ or,
equivalently, $u+v+w\in\mathcal{U}$ then
$\bar{u}+\bar{v}=\overline{u+v}\in\bar{\mathcal{U}}$ iff
$u+v\in\mathcal{U}$. Since $\mathcal{U}$ is ring-supporting, the
latter happens iff $v+w\in\mathcal{U}$ iff
$\overline{v+w}=\bar{v}+\bar{w}\in\bar{\mathcal{U}}$. Therefore
$\bar{\mathcal{U}}$ is a ring-supporting subset of $\mathbf{Z}_n$.

$2)\Longrightarrow 1)$ Go backward in the other implication.

Since $p^{-1}(\bar{\mathcal{U}})=\mathcal{U}$, we immediately get
$\mathcal{U}=\bigcup_{i\in\bar{\mathcal{U}}}(n\mathbf{Z}+i)$, thus
ending the proof.

\end{proof}

\begin{teor} \label{ring-supporting}
Let $\mathcal{U}$ be a  ring-supporting subset of $\mathbf{Z}$
such that $(\mathcal{U}:\mathcal{U})=n\mathbf{Z}$. Then
$\mathcal{U}=\bigcup_{j\in J}(n\mathbf{Z}+j)$, where $J$ is a
ring-supporting subset of $\mathbf{Z}_n$ such that $(J:J)=J\cap
(-J)=\{\bar{0}\}$. Moreover, the
 following
assertions hold:

\begin{enumerate}
\item If $\mathcal{U}$ contains an infinite interval then exactly
one of the following conditions hold:

\begin{enumerate}
\item $\mathcal{U}=\mathbf{Z}=(\mathcal{U}:\mathcal{U})$ \item
$\mathcal{U}$ is an additive submonoid of $\mathbf{N}$ \item
$\mathcal{U}$ is an additive submonoid of $-\mathbf{N}$
\end{enumerate}
 \item
If $\mathcal{U}$ does not contain an infinite interval, then
$\mathcal{U}$ is a union $\mathcal{U}=\bigcup_{i\in I}[a_i,b_i]$
of finite intervals such that $a_i\leq b_i<a_{i+1}-1$ for all
$i\in I$, where the index set $I$ is an interval of $\mathbf{Z}$.
 Moreover, the interval $[a_{i_0},b_{i_0}]$ containing $0$ is either $[0,r]$ or
$[-r,0]$, for some natural number $r$ which, in case $n>1$,
satisfies that $2r<n$
\end{enumerate}

\end{teor}

\begin{proof}

The first part of the theorem follows directly from Lemma
\ref{reduction-modulo-(U:U)}. To prove assertion 1, notice that,
since $(\mathcal{U}:\mathcal{U})$ is a subgroup of $\mathbf{Z}$
contained in $\mathcal{U}$, one has that
$(\mathcal{U}:\mathcal{U})=\mathbf{Z}$ if and only if
$\mathcal{U}=\mathbf{Z}$. So in the rest of the proof we assume
that $\mathcal{U}$ is a proper subset of $\mathbf{Z}$. Suppose
that $\mathcal{U}$  contains an infinite interval $[m,+\propto )$.
  If $m<0$ then the subgroup
$(\mathcal{U}:\mathcal{U})=\mathcal{U}\cap (-\mathcal{U})$
contains the interval $[m,-m]$, which contains $1$. That implies
that $(\mathcal{U}:\mathcal{U})=\mathbf{Z}=\mathcal{U}$, which is
a contradiction. Therefore $m\geq 0$ and we choose $m$ minimal
with the property that $[m,+\propto )\subseteq\mathcal{U}$. Notice
that if $m>0$ then $m>1$, because $0\in\mathcal{U}$.  We also have
$(\mathcal{U}:\mathcal{U})=0$, for otherwise the subgroup
$(\mathcal{U}:\mathcal{U})$ would be lower and upper unbounded and
then $\mathbf{Z}=(\mathcal{U}:\mathcal{U})+[m,+\propto )\subseteq
(\mathcal{U}:\mathcal{U})+\mathcal{U}=\mathcal{U}$, which
contradicts the fact that $\mathcal{U}$ is proper. Then
$\mathcal{U}\cap(-\mathcal{U})=(\mathcal{U}:\mathcal{U})=0$ and,
in particular, $(-\propto ,-m]\cap\mathcal{U}=\emptyset$ unless
$m=0$, in which case that intersection is  $\{0\}$. Then
$\mathcal{U}\subseteq [-m,+\propto )$. In order to prove that
$\mathcal{U}\subset\mathbf{N}$ it will be enough to consider the
case $m>1$ and prove that $(-m ,0)\cap\mathcal{U}=\emptyset$.
Indeed, if  that intersection is nonempty, then we can take $k>0$
such that  $-k=Inf[(-m ,0)\cap\mathcal{U}]$. We choose now
$u=v=-k$ and $w=m+2k$. Then $u,v,w$ and $u+v+w=m$ belong to
$\mathcal{U}$. However $v+w=m+k\in\mathcal{U}$ while
$u+v=-2k\not\in\mathcal{U}$. That contradicts the ring-supporting
condition of $\mathcal{U}$. Therefore
$\mathcal{U}\subseteq\mathbf{N}$. We finally check that
$\mathcal{U}$ is an additive submonoid of $\mathbf{N}$. Suppose
that it is not the case, so that there are $u,v\in\mathcal{U}$
such that $u+v\not\in\mathcal{U}$. Then, taking $w=m$, we get that
$u,v,w\in\mathcal{U}$ and $u+v+w\geq m$, so that also
$u+v+w\in\mathcal{U}$. But now $v+w=m+v\in
[m,+\propto)\subseteq\mathcal{U}$ and $u+v\not\in\mathcal{U}$.
That contradicts the ring-supporting condition of $\mathcal{U}$.

Since $-\mathcal{U}$ is a ring-supporting subset whenever so is
$\mathcal{U}$ (cf. Proposition \ref{ejemplos-modular-pairs}), the
above paragraph also gives that if $(-\propto
,-m]\subseteq\mathcal{U}$, for some $m\in\mathbf{Z}$, then
$\mathcal{U}$ is an additive submonoid of $-\mathbf{N}$.

 Suppose now that  $\mathcal{U}$ does not
contain infinite intervals. We can always express it as a union
$\mathcal{U}=\bigcup_{i\in I}[a_i,b_i]$  of finite intervals such
that $a_i\leq b_i<a_{i+1}-1$ for all $i\in I$, where the index set
$I$ is an interval of the integers which we convene that contains
$0$ and, in addition $0\in [a_0,b_0]$.
 Suppose now that $a_0<0$ and
$b_0>0$. We then take $u=a_0$, $v=-1$, $w=1$. Then $u,v,w\in
[a_0,b_0]\subset\mathcal{U}$, and also $u+v+w=u\in\mathcal{U}$.
Now $v+w=0\in\mathcal{U}$ but $u+v=a_0-1\notin\mathcal{U}$. That
contradicts the ring-supporting condition of $\mathcal{U}$.
Therefore either $a_0=0$ or $b_0=0$, so that we can assume that
$[a_0,b_0]=[0,r]$ or $[a_0,b_0]=[-r,0]$, for some natural number
$r$. It only remains to prove that, in case $n>1$, we have $2r<n$.
We do it when $[a_0,b_0]=[0,r]$, leaving the dualization when
$[a_0,b_0]=[-r,0]$ for the reader. Since $\bar{\mathcal{U}}$ is a
proper subset of $\mathbf{Z}_n\cong [0,n)$, we necessarily have
$r<n-1$. Suppose now that $2r\geq n$. Then we  take $u=n-r$, $v=r$
and $w=1$. Since $0<n-r\leq 2r-r=r$, and $u+v+w=n+1$, we conclude
that $u,v,w$ and $u+v+w$ belong to $\mathcal{U}$. Now
$u+v=n\in\mathcal{U}$, but $v+w=r+1\notin\mathcal{U}$ because
$r=b_0<a_1-1$. That contradicts the fact that $\mathcal{U}$ is
ring-supporting and, hence, we necessarily have $2r<n$.

\end{proof}

\begin{rem}
From the above result we deduce that, in order to classify all the
ring-supporting subsets of $\mathbf{Z}$, one should classify all
the ring-supporting subsets $\mathcal{U}$ such that
$(\mathcal{U}:\mathcal{U})=0$ and, for every $n>0$, all the
ring-supporting subsets $J\subset\mathbf{Z}_n$ such that
$(J:J)=\{\bar{0}\}$.  It seems to be a very difficult task, and
the only reasonable expectation is to tackle the case when $n>0$
is not large. The following is a sample. We leave the verification
to the reader.
\end{rem}

\begin{ejems}
The following are the ring-supporting subsets
$\mathcal{U}\subseteq\mathbf{Z}$ such that
$(\mathcal{U}:\mathcal{U})=n\mathbf{Z}$, with $0<n\leq 5$

\begin{enumerate}
\item If $n=1$ then $\mathcal{U}=\mathbf{Z}$ \item If $n=2$ then
$\mathcal{U}=2\mathbf{Z}=\bigcup_{k\in\mathbf{Z}}[2k,2k]$ \item If
$n=3$ then $\mathcal{U}$ is one of the subsets
$3\mathbf{Z}=\bigcup_{k\in\mathbf{Z}}[3k,3k]$, $3\mathbf{Z}\cup
(3\mathbf{Z}+1)=\bigcup_{k\in\mathbf{Z}}[3k,3k+1]$ or
$3\mathbf{Z}\cup
(3\mathbf{Z}+2)=\bigcup_{k\in\mathbf{Z}}[3k-1,3k]$ \item If $n=4$
then $\mathcal{U}$ is one of the subsets
$4\mathbf{Z}=\bigcup_{k\in\mathbf{Z}}[4k,4k]$, $4\mathbf{Z}\cup
(4\mathbf{Z}+1)=\bigcup_{k\in\mathbf{Z}}[4k,4k+1]$ or
$4\mathbf{Z}\cup
(4\mathbf{Z}+3)=\bigcup_{k\in\mathbf{Z}}[4k-1,4k]$ \item If $n=5$
then $\mathcal{U}$ is one of the subsets of the following list:

\begin{enumerate}
\item[a)]$5\mathbf{Z}=\bigcup_{k\in\mathbf{Z}}[5k,5k]$

\item[b)] $5\mathbf{Z}\cup
(5\mathbf{Z}+1)=\bigcup_{k\in\mathbf{Z}}[5k,5k+1]$ \item[b*)]
$5\mathbf{Z}\cup
(5\mathbf{Z}+4)=\bigcup_{k\in\mathbf{Z}}[5k-1,5k]$ \item[c)]
 $5\mathbf{Z}\cup
(5\mathbf{Z}+1)\cup
(5\mathbf{Z}+2)=\bigcup_{k\in\mathbf{Z}}[5k,5k+2]$ \item[c*)]
$5\mathbf{Z}\cup (5\mathbf{Z}+3)\cup
(5\mathbf{Z}+4)=\bigcup_{k\in\mathbf{Z}}[5k-2,5k]$ \item[d)]
$5\mathbf{Z}\cup
(5\mathbf{Z}+2)=\bigcup_{k\in\mathbf{Z}}([5k,5k]\cup [5k+2,5k+2])$
\item [d*)] $5\mathbf{Z}\cup
(5\mathbf{Z}+3)=\bigcup_{k\in\mathbf{Z}}([5k,5k]\cup [5k+3,5k+3])$
\item[e)] $5\mathbf{Z}\cup (5\mathbf{Z}+1)\cup
(5\mathbf{Z}+3)=\bigcup_{k\in\mathbf{Z}}([5k,5k+1]\cup
[5k+3,5k+3])$ \item[e*)] $5\mathbf{Z}\cup (5\mathbf{Z}+2)\cup
(5\mathbf{Z}+4)=\bigcup_{k\in\mathbf{Z}}([5k-1,5k]\cup
[5k+2,5k+2])$
\end{enumerate}

\end{enumerate}
\end{ejems}

The last part of  Theorem \ref{ring-supporting} gives the
following  peculiar ring-supporting subset, which will be called
{\bf translation of an interval} in the sequel:

\begin{cor} \label{translation-of-interval}
Let $\mathcal{U}$ be a  subset of $\mathbf{Z}$ containing $0$ such
that $(\mathcal{U}:\mathcal{U})=n\mathbf{Z}$, with $n\neq 0,1$.
The following assertions are equivalent:

\begin{enumerate}
\item $\mathcal{U}$ is ring-supporting and
$\mathcal{U}=I+n\mathbf{Z}$, for some finite interval
$I\subset\mathbf{Z}$ \item
$\mathcal{U}=\bigcup_{k\in\mathbf{Z}}[nk,nk+r]$ or
$\mathcal{U}=\bigcup_{k\in\mathbf{Z}}[nk-r,nk]$, for some natural
number $r$ such that $2r<n$
\end{enumerate}
\end{cor}

\begin{proof}
Theorem \ref{ring-supporting} gives the implication
$1)\Longrightarrow 2)$. For the converse, it will enough to check
that $\mathcal{U}=\bigcup_{k\in\mathbf{Z}}[nk,nk+r]$ is
ring-supporting and, for that, we will check that
$\bar{\mathcal{U}}$ is a ring-supporting subset of $\mathbf{Z}_n$.
Indeed, if $u,v,w$ are natural numbers such that $u,v,w\leq r$ and
we assume that $\bar{u}+\bar{v}+\bar{w}\in\bar{\mathcal{U}}$,
then, since $u+v\leq 2r<n$, we have that either $u+v+w\leq r$ or
$n\leq u+v+w\leq n+r$. In the first case, both $u+v\leq r$ and
$v+w\leq r$. In the second case, both $u+v>r$ and $v+w>r$ for if,
say, $u+v\leq r$ then $u+v+w\leq r+r=2r<n$, against the assumption
that $u+v+w\in [n,n+r]$. Therefore
$\bar{u}+\bar{v}\in\bar{\mathcal{U}}$ if, and only if,
$\bar{v}+\bar{w}\in\bar{\mathcal{U}}$, as desired.
\end{proof}

In the rest of the section, we assume that our algebra is a
positively graded algebra $A=\oplus_{i\geq 0}A_i$ over a field
$K$,  which is generated in degrees $0,1$ and has the property
that $A_0$ is semisimple. The symbol $\otimes$ will  always mean
tensor product over $A_0$.

In case $(\mathcal{S},\mathcal{U})$ is a right modular pair of
subsets of $\mathbf{Z}$ and $X\in Gr_{A_\mathcal{U}}$ has
$Supp(X)\subseteq\mathcal{S}$, we denote by $\mu_{s,u}:X_s\otimes
A_u\longrightarrow X_{s+u}$ the multiplication map, whenever
$(s,u)\in\mathcal{S}\times\mathcal{U}$ and $s+u\in\mathcal{S}$.
Since $X\otimes A$ is a graded right $A$-module,
$Ker(\mu_{s,u})A_i$ is a well-defined $A_0$-submodule of
$X_s\otimes A_{u+i}$, for all $i\geq 0$. With this terminology in
mind, we can give criteria for $X$ to be liftable.

\begin{prop} \label{liftable}
Let $A=\oplus_{i\geq 0}A_i$ be a positively graded algebra
generated in degrees $0,1$, with $A_0$ semisimple,
$(\mathcal{S},\mathcal{U})$ be a right modular pair of subsets of
$\mathbf{Z}$  and consider an $X\in Gr_{A_\mathcal{U}}$ generated
in degrees belonging to $(\mathcal{S}:\mathcal{U})$.
 The following assertions
are equivalent:

\begin{enumerate}
\item $X$ is liftable with respect to the functor
$(-)_\mathcal{S}:Gr_A\longrightarrow Gr_{A_\mathcal{U}}$ \item
$Ker(\mu_{m,u})A_{v-u}\subseteq Ker(\mu_{m,v})$ whenever
$(m,u,v)\in
(\mathcal{S}:\mathcal{U})\times\mathcal{U}\times\mathcal{U}$ is a
triple such that $0\leq u<v$ and $v-u\not\in\mathcal{U}$
\end{enumerate}
\end{prop}

\begin{proof}
$1)\Longrightarrow 2)$ Take $M\in Gr_A$ such that
$M_\mathcal{S}=X$. Then, for all $(m,u)\in
(\mathcal{S}:\mathcal{U})\times\mathcal{U}$,  the multiplication
maps $\mu_{m,u}:X_m\otimes A_u=M_m\otimes A_u\longrightarrow
M_{m+u}=X_{m+u}$ are those of the $A$-module $M$. Now condition 2
follows from the associative property of the multiplication for
the $A$-module $M$

$2)\Longrightarrow 1)$ When $(\mathcal{S}:\mathcal{U})=\emptyset$
there is nothing to prove. Hence, we assume that
$(\mathcal{S}:\mathcal{U})\neq\emptyset$. Since $X$ is generated
in degrees belonging to $(\mathcal{S}:\mathcal{U})$, the
multiplication map $\mu :X_{(\mathcal{S}:\mathcal{U})}\otimes
A_\mathcal{U}\longrightarrow X$ is surjective and
$Supp(X)\subseteq\mathcal{S}$. We then consider the graded
$A$-module $M=\frac{X_{(\mathcal{S}:\mathcal{U})}\otimes A}{(Ker
\mu )A}$. Then
$M_\mathcal{S}=\frac{(X_{(\mathcal{S}:\mathcal{U})}\otimes
A)_\mathcal{S}}{[(Ker \mu
)A]_\mathcal{S}}=\frac{X_{(\mathcal{S}:\mathcal{U})}\otimes
A_\mathcal{U}}{[(Ker \mu )A]_\mathcal{S}}$. We shall prove that
$[(Ker \mu)A]_\mathcal{S}=Ker\mu$ and the result will follow, for
$\frac{X_{(\mathcal{S}:\mathcal{U})}\otimes A_\mathcal{U}}{Ker \mu
}$ is isomorphic to $X$ in $Gr_{A_\mathcal{U}}$. Let us fix
$t\in\mathcal{S}$. Then the $t$-homogeneous component of $(Ker\mu
)A$ is $\sum_{s\in\mathcal{S}, s\leq t}(Ker\mu )_sA_{t-s}$. The
task is whence reduced to prove that if $s\leq t$ are elements of
$\mathcal{S}$ then $(Ker\mu )_sA_{t-s}\subseteq (Ker\mu )_t$. That
inclusion is clear in case $t-s\in\mathcal{U}$ due to the
associative property of the multiplication in the
$A_\mathcal{U}$-module  $X$. So we assume that
$t-s\not\in\mathcal{U}$. Let $\sum_{1\leq i\leq r}x_i\otimes a_i$
belong to $(Ker\mu )_s$, where $deg(x_i)=m_i\in
(\mathcal{S}:\mathcal{U})$ and $a_i\in A_{s-m_i}$ for $i=1,...,r$.
We take $p=max\{m\in (\mathcal{S}:\mathcal{U}):$ $m\leq s\}$. The
hypothesis on $A$ implies that $A_{s-m_i}=A_{p-m_i}\cdot A_{s-p}$,
so that $a_i=\sum_jb_{ij}c_{ij}$, where $b_{ij}\in A_{p-m_i}$ and
$c_{ij}\in A_{s-p}$, for all $i,j$. Then
$0=\sum_ix_ia_i=\sum_{i,j}x_i(b_{ij}c_{ij})=\sum_{i,j}(x_ib_{ij})c_{ij}$,
using the associative property of $X$ again (notice that
$deg(b_{ij})=p-m_i$ and $deg(c_{ij})=s-p$ belong to
$\mathcal{U}$). Then $\sum_{i,j}(x_ib_{ij})\otimes c_{ij}$ belongs
to the kernel of the multiplication map $\mu_{p,s-p}:X_p\otimes
A_{s-p}\longrightarrow A_s$. Since $p\in
(\mathcal{S}:\mathcal{U})$, we have that $u=s-p$ and $v=t-p$ are
elements of $\mathcal{U}$ such that $0\leq u<v$ and
$v-u=t-s\not\in\mathcal{U}$. Our hypothesis says that
$Ker(\mu_{p,s-p})A_{t-s}\subseteq Ker(\mu_{p,t-p})$. Then, if
$d\in A_{t-s}$ the element $\sum_{i,j}(x_ib_{ij})\otimes
(c_{ij}d)$ belongs to $Ker(\mu_{p,t-p})$. That gives
$0=\sum_{i,j}(x_ib_{ij})(c_{ij}d)=\sum_{i,j}x_{i}[b_{ij}(c_{ij}d)]$,
using the associative property of the multiplication in $X$. But
the associative property of $A$ gives that
$b_{ij}(c_{ij}d)=(b_{ij}c_{ij})d$. Hence we have
$0=\sum_{i,j}x_i[(b_{ij}c_{ij})d]=\sum_ix_i(a_id)$, which implies
that the element $(\sum_ix_i\otimes a_i)d$ belongs to
$(Ker\mu)_t$. That proves that $(Ker\mu_s)A_{t-s}\subseteq (Ker\mu
)_t$, as desired.

\end{proof}

 The conditions get much simpler in the case that we
are more interested in, namely, the translation of an interval
(cf. Corollary \ref{translation-of-interval}).

\begin{cor} \label{liftable-interval}
In the situation of the proposition, adopting the notations of
Corollary \ref{translation-of-interval}, suppose that
$\mathcal{U}$ is the translation of an interval. Then the
following assertions hold for a graded $A_\mathcal{U}$-module $X$
generated in degrees belonging to $(\mathcal{S}:\mathcal{U})$:

\begin{enumerate}
\item If $\mathcal{U}=\bigcup_{k\in\mathbf{Z}}[kn,kn+r]$ ($0\leq
2r<n$), then $X$ is liftable if, and only if,
$Ker(\mu_{m,r})A_{n-r}\subseteq Ker(\mu_{m,n})$, for all $m\in
(\mathcal{S}:\mathcal{U})$ \item If
$\mathcal{U}=\bigcup_{k\in\mathbf{Z}}[kn-r,kn]$ ($0\leq 2r<n$),
then $X$ is liftable if, and only if,
$Ker(\mu_{m,u})A_{v-u}\subseteq Ker(\mu_{m,v})$, for all $m\in
(\mathcal{S}:\mathcal{U})$ and all $n-r\leq u<v\leq n$
\end{enumerate}
\end{cor}

\begin{proof}
We prove 1, for which we only need to check the 'if' part. The
proof of assertion 2 follows a similar pattern. We first prove by
induction on $k\geq 0$ that $Ker(\mu_{m,kn+r})A_{n-r}\subseteq
Ker(\mu_{m,(k+1)n})$ (*), for all $m\in
(\mathcal{S}:\mathcal{U})$. The case $k=0$ is just the hypothesis.
Suppose  now  that $k>0$ and $\sum_ix_i\otimes a_i\in
Ker(\mu_{m,kn+r})$, where $x_i\in X_m$ and $a_i\in A_{kn+r}$ for
all indices $i$. Since $A_{kn+r}=A_n\cdot A_{(k-1)n+r}$, we can
decompose $a_i=\sum_jb_{ij}c_{ij}$, where $b_{ij}\in A_n$ and
$c_{ij}\in A_{(k-1)n+r}$. With an argument similar to that in the
proof of Proposition \ref{liftable}, we see that
$\sum_{i,j}(x_ib_{ij})\otimes c_{ij}$ belongs to
$Ker(\mu_{m+n,(k-1)n+r})$. If now $d\in A_{n-r}$ then the
induction hypothesis says that $(\sum_{i,j}(x_ib_{ij})\otimes
c_{ij})d=\sum_{i,j}(x_ib_{ij})\otimes (c_{ij}d)$ belongs to
$Ker(\mu_{m+n,kn})$ and, as in the proof of Proposition
\ref{liftable} again, one sees that $(\sum_ix_i\otimes
a_i)d=\sum_i x_i\otimes (a_id)\in Ker(\mu_{m,(k+1)n})$.

Let $(m,u,v)\in
(\mathcal{S}:\mathcal{U})\times\mathcal{U}\times\mathcal{U}$ be a
triple such that $0\leq u<v$ and $v-u\not\in\mathcal{U}$. Then we
have $u\in [a_i,b_i]=[in,in+r]$ and $v\in [a_j,b_j]=[jn,jn+r]$,
for some indices $0\leq i<j$. Since $A$ is generated in degrees
$0,1$, we have that $A_{v-u}=A_{b_i-u}\cdot A_{a_{i+1}-b_i}\cdot
...\cdot A_{a_j-b_{j-1}}\cdot A_{v-a_j}$. Each of the factors of
this product is of the form $A_{v'-u'}$, where either $u',v'\in
[a_k,b_k]=[kn,kn+r]$ or $u'=b_k=kn+r$, $v'=a_{k+1}=(k+1)n$, for
some index $k$. In the first case $0\leq v'-u'\leq r$ and, hence,
$v'-u'\in\mathcal{U}$. In the second case, $v'-u'=n-r$. One then
obtains that $Ker(\mu_{m,u})A_{v-u}\subseteq Ker(\mu_{m,v})$ by
recursive use of the associative property of the
$A_\mathcal{U}$-module $X$ and the above condition (*).
\end{proof}

Using Theorem \ref{fullyfaithful-killing-supports}, we have the
following straightforward consequence:

\begin{cor} \label{equivalence-for-translation-of-intervals}
Let $A=\oplus_{i\geq 0}A_i$ be a positively graded algebra
generated in degrees $0,1$ such that $A_0$ is semisimple. We
consider  $\mathcal{U}=\bigcup_{k\in\mathbf{Z}}[nk,nk+r]$ and
$\mathcal{S}=m+\mathcal{U}$, where $n>1$,  $0\leq 2r<n$ and
$m\in\mathbf{Z}$. The functor $(-)_\mathcal{S}:Gr_A\longrightarrow
Gr_{A_\mathcal{U}}$ establishes an equivalence  between the
following categories:

\begin{enumerate}
\item The full subcategory $\mathcal{G}(\mathcal{S},\mathcal{U})$
of $Gr_A$ whose objects are the graded modules generated in
degrees belonging to $m+n\mathbf{Z}$ and cogenerated in degrees
belonging to $\mathcal{S}$ \item The full subcategory
$\mathcal{L}(\mathcal{S},\mathcal{U})$ of $Gr_{A_\mathcal{U}}$
whose objects are those $X$ generated in degrees belonging to
$m+n\mathbf{Z}$ and satisfying that
$(Ker\mu_{m+kn,r})A_{n-r}\subseteq Ker\mu_{m+kn,n}$  for all
$k\in\mathbf{Z}$
\end{enumerate}
\end{cor}

By the induction principle, given a subset
$\mathcal{S}\subseteq\mathbf{Z}$, which is neiter upper nor lower
bounded,  and a fixed element $s_0\in\mathcal{S}$, there is a
unique strictly increasing map $\delta
=\delta_{(\mathcal{S},s_0)}:\mathbf{Z}\longrightarrow\mathbf{Z}$
such that $\delta (0)=s_0$ and $Im(\delta )=\mathcal{S}$.

Let $(\mathcal{S},\mathcal{U})$ be a right modular pair of subsets
  of $\mathbf{Z}$ such that $0\neq
(\mathcal{U}:\mathcal{U})=n\mathbf{Z}$ and $m\in
(\mathcal{S}:\mathcal{U})$. We want to identify
$\delta_{(\mathcal{S},m)}$. As in Lemma
\ref{reduction-modulo-(U:U)}, we can think of $\bar{\mathcal{U}}$
as a subset of $[0,n)$, in which case we can write
$\bar{\mathcal{U}}=\{i_0,i_1,...,i_{t-1}\}$, with
$t=|\bar{\mathcal{U}}|$ and $0=i_0<i_1<...<i_{t-1}<n$.

\begin{lema} \label{funcion asociada a un ring-supporting}
In the above situation, the following assertions hold:

\begin{enumerate}

\item $\delta_{(\mathcal{S},m)}=m+\delta_{(\mathcal{U},0)}$ \item
$\delta_{(\mathcal{S},m)}:\mathbf{Z}\longrightarrow\mathbf{Z}$
maps $tj+k\rightsquigarrow m+nj+i_k$, for all $j\in\mathbf{Z}$ and
all $k\in [0,t)$ \item If $\mathcal{U}$ is the translation of an
interval, then  $\delta_{(\mathcal{U},0)}$ is a pseudomorphism of
groups.
\end{enumerate}
\end{lema}

\begin{proof}

1) Clearly, $m+\delta_{(\mathcal{U},0)}:j\rightsquigarrow
m+\delta_{(\mathcal{U},0)}(j)$ is strictly increasing and maps
$0\rightsquigarrow m+\delta_{(\mathcal{U},0)}(0)=m$. Moreover,
$Im(m+\delta_{(\mathcal{U},0)})=m+Im(\delta_{(\mathcal{U},0)})=m+\mathcal{U}=\mathcal{S}$.
Then $\delta_{(\mathcal{S},m)}=m+\delta_{(\mathcal{U},0)}$.

2) Use an argument as in 1.

3)  We then prove that $\delta=:\delta_{(\mathcal{U},0)}$ is a
pseudomorphism in case
$\mathcal{U}=\bigcup_{k\in\mathbf{Z}}[nk,nk+r]$ or
$\mathcal{U}=\bigcup_{k\in\mathbf{Z}}[nk-r,nk]$, where $r$ and $n$
are  natural numbers such that $0\leq 2r<n$. We consider the
second possibility, leaving the first one as an exercise. Notice
that $t=|\bar{\mathcal{U}}|=r+1$ in this case. Then,  according to
assertion 2, we have $\delta [(r+1)j-i]=nj-i$, for all
$j\in\mathbf{Z}$ and $i\in [0,r]$. If $m=(r+1)j-i$ and
$m'=(r+1)j'-i'$, where $i,i'\in [0,r]$, then $\delta (m)+\delta
(m')=n(j+j')-i-i'$ belongs to $\mathcal{U}=Im(\delta )$ if, and
only if, $i+i'\in [0,r]$. In that case, $\delta (m+m')=\delta
[(r+1)(j+j')-(i+i')]=n(j+j')-(i+i')=\delta (m)+\delta (m')$. Then
$\delta$ is a pseudomorphism of groups.
\end{proof}

We can now combine the general results of Sections 2 and 3 with
those of this section. If
$\mathcal{U}=\bigcup_{k\in\mathbf{Z}}[nk,nk+r]$, with $0<2r<n$,
then we consider the graded algebra $\tilde{B}$ obtained from
$B=:A_\mathcal{U}$ by regrading along the pseudomorphism $\delta
=\delta_{(\mathcal{U},0)}$. If now $V\in Gr_{\tilde{B}}$, then we
have a multiplication maps
$\tilde{\mu}_{(r+1)k,r}:V_{(r+1)k}\otimes\tilde{B}_r=V_{(r+1)k}\otimes
A_r\longrightarrow V_{(r+1)(k+1)-1}$ and
$\tilde{\mu}_{(r+1)k,r+1}:V_{(r+1)k}\otimes\tilde{B}_{r+1}=V_{(r+1)k}\otimes
A_n\longrightarrow V_{(r+1)(k+1)}$, so that, working in the right
$A$-module $V_{(r+1)\mathbf{Z}}\otimes A$, it makes sense to
consider $Ker(\tilde{\mu})_{(r+1)k,r}A_{n-r}$, which will be
$A_0-A_0-$sub-bimodule of
$V_{(r+1)k}\otimes\tilde{B}_{r+1}=V_{(r+1)k}\otimes A_n$.  We then
get the following:

\begin{cor}
Let $A=\oplus_{n\geq 0}A_n$ be a positively graded algebra
generated in degrees $0,1$ such that $A_0$ is semisimple,
$(\mathcal{S}, \mathcal{U})=(m+\mathcal{U},\mathcal{U})$ be a
right modular pair, where
$\mathcal{U}=\bigcup_{k\in\mathbf{Z}}[nk,nk+r]$, with $0\leq
2r<n$, and $\tilde{B}$ be the graded algebra obtained from
$B=A_\mathcal{U}$ by regrading along the pseudomorphism $\delta
=\delta_{(\mathcal{U},0)}:\mathbf{Z}\longrightarrow\mathbf{Z}$.
The functor $\Phi :Gr_A\longrightarrow Gr_{\tilde{B}}$ of
Proposition \ref{ring-supporting=image of pseudomorphism} gives an
equivalence of categories between:

\begin{enumerate}
\item The full subcategory $\mathcal{G}(\mathcal{S},\mathcal{U})$
of $Gr_A$ whose objects are the graded $A$-modules generated in
degrees belonging to $(\mathcal{S}:\mathcal{U})=m+n\mathbf{Z}$ and
cogenerated in degrees belonging to $\mathcal{S}$ \item The full
subcategory $\mathcal{L}_{\tilde{B}}$ of $Gr_{\tilde{B}}$
consisting of those $V$ generated in degrees belonging to
$(r+1)\mathbf{Z}$ such that
$Ker(\tilde{\mu}_{(r+1)k,r})A_{n-r}\subseteq
Ker(\tilde{\mu}_{(r+1)k,r+1})$, for all $k\in\mathbf{Z}$.
\end{enumerate}

Moreover, $\mathcal{L}_{\tilde{B}}$ contains all the graded
$\tilde{B}$-modules presented in degrees belonging to
$(r+1)\mathbf{Z}$.
\end{cor}

\begin{proof}
It is a direct consequence of Proposition
\ref{ring-supporting=image of pseudomorphism} and Corollary
\ref{equivalence-for-translation-of-intervals}, after checking
that the equivalence $\Sigma_\mathcal{S}=\Sigma_{m+Im(\delta
)}\stackrel{\cong}{\longrightarrow}\tilde{\Sigma}_{Im(\delta
)}=\tilde{\Sigma}_\mathcal{U}$ of Proposition
\ref{caracterizacion-de-pseudomorphisms} restricts to an
equivalence between $\mathcal{L}(\mathcal{S},\mathcal{U})$ and
$\mathcal{L}_{\tilde{B}}$.
\end{proof}

\begin{rem} \label{n-Koszul}
If $\Lambda$ is a $n$-Koszul algebra, with $n>2$, then we can
apply the above corollary to $A=\Lambda^!$, with $r=1$, so that
$\tilde{B}\cong E=E(\Lambda )$ is the Yoneda algebra of $\Lambda$
(cf. \cite{GMMVZ}[Theorem 9.1], see also \cite{BM}[Proposition
3.1]) and $(\mathcal{S},\mathcal{U})=(m+\mathcal{U},\mathcal{U})$,
with $\mathcal{U}=n\mathbf{Z}\cup (n\mathbf{Z}+1)$. Hence, if
$V\in Gr_E$, we have multiplication maps
$\tilde{\mu}_{2k,1}:V_{2k}\otimes
E_1=V_{2k}\otimes\Lambda^!_1\longrightarrow V_{2k+1}$ and
$\tilde{\mu}_{2k,2}:V_{2k}\otimes
E_2=V_{2k}\otimes\Lambda^!_n\longrightarrow V_{2k+2}$ and we can
take $Ker(\tilde{\mu}_{2k,1})\Lambda^!_{n-1}$, which is a
$\Lambda_0-\Lambda_0-$sub-bimodule of
$V_{2k}\otimes\Lambda^!_n=V_{2k}\otimes E_2$. The last corollary
gives then an equivalence of categories between:

\begin{enumerate}
\item The full subcategory $\mathcal{G}(\mathcal{S},\mathcal{U})$
of $Gr_{\Lambda^!}$ consisting of those $M$ such that
$M=M_{m+n\mathbf{Z}}\Lambda^!$ and, for every graded submodule
$0\neq N<M$, one has $N_\mathcal{S}\neq 0$ \item The full
subcategory $\mathcal{L}_E$ of $Gr_E$ consisting of those $V\in
Gr_E$ which are generated in even degrees and satisfy that
$Ker(\tilde{\mu}_{2k,1})\Lambda^!_{n-1}\subseteq
Ker(\tilde{\mu}_{2k,2})$, for all $k\in\mathbf{Z}$.
\end{enumerate}

Moreover, $\mathcal{L}_E$ contains all the graded $E$-modules
presented in even degrees.
\end{rem}

\end{document}